\documentstyle[12pt]{article}
\begin{document}
\begin{titlepage}
\begin{flushright}
UQMATH-qe2-9701\\
q-alg/9703020
\end{flushright}
\vskip.3in

\begin{center}
{\Large \bf Comments on Drinfeld Realization of Quantum Affine
Superalgebra $U_q[gl(m|n)^{(1)}]$ and Its Hopf Algebra Structure}
\vskip.3in
{\large Yao-Zhong Zhang} \footnote{Queen Elizabeth II Fellow;
Email: yzz@maths.uq.edu.au}
\vskip.2in
{\em Department of Mathematics, University of Queensland, Brisbane,
Qld 4072, Australia}
\end{center}

\vskip 2cm
\begin{center}
{\bf Abstract}
\end{center}
By generalizing the Reshetikhin and Semenov-Tian-Shansky construction
to supersymmetric cases, we obtain Drinfeld current realization
for quantum affine superalgebra $U_q[gl(m|n)^{(1)}]$. We find a simple
coproduct for the quantum current generators and establish the Hopf
algebra structure of this super current algebra.


\end{titlepage}


\def\a{\alpha}
\def\b{\beta}
\def\d{\delta}
\def\e{\epsilon}
\def\g{\gamma}
\def\k{\kappa}
\def\l{\lambda}
\def\o{\omega}
\def\t{\theta}
\def\s{\sigma}
\def\D{\Delta}
\def\L{\Lambda}

\def\G{{\cal G}}
\def\C{{\bf C}}
\def\P{{\bf P}}

\def\uqgh{{U_q[gl(m|n)^{(1)}]}}
\def\uqg{{U_q[gl(m|n)]}}


\def\beq{\begin{equation}}
\def\eeq{\end{equation}}
\def\bea{\begin{eqnarray}}
\def\eea{\end{eqnarray}}
\def\ba{\begin{array}}
\def\ea{\end{array}}
\def\no{\nonumber}
\def\lt{\left}
\def\rt{\right}
\newcommand{\bq}{\begin{quote}}
\newcommand{\eq}{\end{quote}}

\newtheorem{Theorem}{Theorem}
\newtheorem{Definition}{Definition}
\newtheorem{Proposition}{Proposition}
\newtheorem{Lemma}[Theorem]{Lemma}
\newtheorem{Corollary}[Theorem]{Corollary}
\newcommand{\proof}[1]{{\bf Proof. }
        #1\begin{flushright}$\Box$\end{flushright}}

\newcommand{\sect}[1]{\setcounter{equation}{0}\section{#1}}
\renewcommand{\theequation}{\thesection.\arabic{equation}}

\sect{Introduction\label{intro}}

This paper contains two results on the Drinfeld second realization 
(current realization) \cite{Dri88}
for the untwisted quantum affine superalgebra $\uqgh$.

The first result is the extension of Reshetikhin and Semenov-Tian-Shansky
(RS) construction \cite{Res90} to supersymmetric cases. Using this
super RS construction and a super version of the Ding-Frenkel 
theorem \cite{Din93}, we obtain the defining relations for $\uqgh$
in terms of super (or graded) current generators (generating functions). 
In \cite{Fan97}, the
authors made a similar effort. However, the relations they obtained
are {\em not}
supersymmetric and the algebra they defined is not a current
realization of $\uqgh$ but rather a current
realization of a ``nonstandard" quantum {\em bosonic} algebra associated with
$U_q[gl(N=m+n)]$. This is because those authors 
failed to take care of the {\em grading} in the multiplication rule of tensor
products which plays a fundamental role in any supersymmetric theories.

The second result is the coproduct, counit and antipode for the current
realization of $U_q[gl(m|n)^{(1)}]$, thus establishing a Hopf algebra 
structure of this super current algebra.

\sect{Super RS Algebra and Ding-Frenkel Theorem}

Let us start with introducing some useful notations.
The graded Yang-Baxter equation (YBE) with spectral-parameter dependence takes 
the form
\beq\label{rrr}
R_{12}({z\over w})R_{13}(z)R_{23}(w)=R_{23}(w)R_{13}(z)R_{12}({z\over w}),
\eeq
where $R(z)\in End(V\otimes V)$ with $V$ being graded vector space
and  obeys the weight conservation condition: 
$R(z)^{\a'\b'}_{\a\b}\neq 0$  only when
$[\a']+[\b']+[\a]+[\b]=0$ mod$2$.
The multiplication rule for the tensor product is defined for
homogeneous elements $a,~b,~c,~d$ of a quantum superalgebra by
\beq
(a\otimes b)(c\otimes d)=(-1)^{[b][c]}\,(ac\otimes bd)
\eeq
where $[a]\in {\bf Z}_2$ denotes the grading of the element $a$. 

Introduce the graded permutation operator $P$ on the tensor product
module $V\otimes V$ such that $P(v_\a\otimes v_\b)=(-1)^{[\a][\b]}
(v_\b\otimes v_\a)\,,~\forall v_\a, v_\b\in V$. 
In most cases R-matrix enjoys, among others, the following properties 
\bea
&& (i)~ P_{12}R_{12}(z)P_{12}=R_{21}(z),\label{pt-symmetry}\\
&& (ii)~ R_{12}({z\over w})R_{21}({w\over z})=1.\label{unitarity}
\eea

The graded YBE, when written in matrix form, 
carries extra signs,  and depending on definition of matrix elements,
these sign factors differ \cite{Zha96}. If we define the matrix elements of
$R(z)$ by
\beq
R(z)(v^{\a'}\otimes v^{\b'})
  =R(z)^{\a'\b'}_{\a\b}(v^\a\otimes v^\b),
\eeq
then the matrix YBE reads \cite{Bra94}
\bea
&&R({z\over w})_{\a\b}^{\a'\b'}
  R(z)_{\a'\g}^{\a''\g'}
  R(w)_{\b'\g'}^{\b''\g''}
  (-1)^{[\a][\b]+[\g][\a']+[\g'][\b']}\nonumber\\
&&~~~~~~=R(w)_{\b\g}^{\b'\g'}
  R(z)_{\a\g'}^{\a'\g''}
  R({z\over w})_{\a'\b'}^{\a''\b''}
  (-1)^{[\b][\g]+[\g'][\a]+[\b'][\a']}.\label{ybe-def2'}
\eea
After the redefinition
\beq
\tilde{R}(z)^{\a'\b'}_{\a\b}
 =R(z)^{\a'\b'}_{\a\b}\,(-1)^{[\a][\b]}\label{redef2}
\eeq
the extra signs in (\ref{ybe-def2'}) disappear. However, this redefinition
does not preserve semiclassical properties.

In matrix form the graded permutation operator
$P=\sum_{\a,\b}\,(-1)^{[\b]}\;E^\a_\b\otimes E^\b_\a$ reads
\beq
P^{\a'\b'}_{\a\b}=\d_{\a\b'}\d_{\a'\b}
  (-1)^{[\a'][\b']}.
\eeq

The RS construction \cite{Res90} can be generalized to supersymmetric cases.
Formally, the super RS algebra is defined by similar relations as in 
non-supersymmetric cases \cite{Res90,Din93}, but tensor products now
carry gradings. We are thus led to

\begin{Definition}\label{rs}:  Super RS algebra is
generated by invertible $L^\pm(z)$, satisfying
\bea
R({z\over w})L_1^\pm(z)L_2^\pm(w)&=&L_2^\pm(w)L_1^\pm(z)R({z\over
         w}),\no\\
R({z_+\over w_-})L_1^+(z)L_2^-(w)&=&L_2^-(w)L_1^+(z)R({z_-\over
         w_+}),\label{super-rs}
\eea
where $L_1^\pm(z)=L^\pm(z)\otimes 1$, $L_2^\pm(z)=1\otimes L^\pm(z)$
and $z_\pm=zq^{\pm {c\over 2}}$. For the first formula of
(\ref{super-rs}), the expansion direction of $R({z\over w})$ can be
chosen in $z\over w$ or $w\over z$, but for the second formula, the
expansion direction must only be in $z\over w$.
\end{Definition}

The super RS algebra is a Hopf algebra: its coproduct is defined by
\beq
\D(L^\pm(z)=L^\pm(zq^{\pm 1\otimes {c\over 2}})\stackrel{.}{\otimes}
     L^\pm(zq^{\mp {c\over 2}\otimes 1}),
\eeq
and its antipode is
\beq
S(L^\pm(z))=L^\pm(z)^{-1}.
\eeq

In matrix form, (\ref{super-rs}) carries extra signs due to the graded
multiplication rule of tensor products:
\bea
&&R({z\over w})_{\a\b}^{\a"\b"}L^\pm(z)_{\a"}^{\a'}L^\pm(w)_{\b"}
      ^{\b'}\,(-1)^{[\a']([\b']+[\b"])}\no\\
&&~~~~~~~~~~~~~~~=L^\pm(w)_\b^{\b"}L^\pm(z)_\a^{\a"}R({z\over w})
      _{\a"\b"}^{\a'\b'}\,(-1)^{[\a]([\b]+[\b"])},\no\\
&&R({z_+\over w_-})_{\a\b}^{\a"\b"}L^+(z)_{\a"}^{\a'}L^-(w)_{\b"}
      ^{\b'}\,(-1)^{[\a']([\b']+[\b"])}\no\\
&&~~~~~~~~~~~~~~~=L^-(w)_\b^{\b"}L^+(z)_\a^{\a"}R({z_-\over w_+})
      _{\a"\b"}^{\a'\b'}\,(-1)^{[\a]([\b]+[\b"])}.\label{rll-component}
\eea
Introduce matrix $\t$:
\beq
\t^{\a'\b'}_{\a\b}=(-1)^{[\a][\b]}\d_{\a\a'}\d_{\b\b'}.\label{theta}
\eeq
With the help of this matrix $\t$, one can cast (\ref{rll-component}) into the
usual matrix equation,
\bea
R({z\over w})L_1^\pm(z)\t L_2^\pm(w)\t&=&\t L_2^\pm(w)\t L_1^\pm(z)R({z\over
         w}),\no\\
R({z_+\over w_-})L_1^+(z)\t L_2^-(w)\t &=&\t L_2^-(w)\t L_1^+(z)R({z_-\over
         w_+}).\label{RLL-LLR1}
\eea
Now the multiplications in ({\ref{RLL-LLR1}) are simply the usual matrix
multiplications.

In this paper we take $R({z\over w})\in End(V\otimes V)$ to be the 
R-matrix associated with
$\uqg$, where $V$ is a $(m+n)$-dimensional graded vector space. Let
basis vectors $\{v^1,\;v^2,\cdots,\;v^m\}$ be even and $\{v^{m+1},\;v^{m+2},
\cdots,\; v^{m+n}\}$ be odd. Then the R-matrix has the following matrix 
elements:
\bea
R({z\over w})_{\a\b}^{\a'\b'}&=&(-1)^{[\a][\b]}\tilde{R}({z\over w})_{\a\b}
  ^{\a'\b'},\no\\
\tilde{R}({z\over w})&=&\sum^m_{i=1}E^i_i\otimes E^i_i
  +\sum_{i=m+1}^{m+n}\frac{wq-zq^{-1}}{zq-wq^{-1}}E^i_i\otimes E^i_i
  +\frac{z-w}{zq-wq^{-1}}\sum_{i\neq j}(-1)^{[i][j]}E^i_i\otimes E^j_j\no\\
& &\sum_{i<j} \frac{z(q-q^{-1})}{zq-wq^{-1}}E^j_i\otimes E^i_j
  +\sum_{i>j}\frac{w(q-q^{-1})}{zq-wq^{-1}}E^j_i\otimes E^i_j.\label{r-matrix}
\eea
It is easy to check that the R-matrix $R(z)$ satisfies
(\ref{pt-symmetry}) and (\ref{unitarity}).
We will construct Drinfeld realization of $\uqgh$.  We first state a super
version of the Ding-Frenkel theorem \cite{Din93}:

\begin{Theorem}\label{df}: $L^\pm(z)$ has the following unique
Gauss decomposition
\bea
L^\pm(z)&=&\left (
\begin{array}{cccc}
1 & \cdots & {} & 0\\
e^\pm_{2,1}(z) & \ddots & {} & {}\\
e^\pm_{3,1}(z) & {}     & {} & \vdots\\
\vdots &  {} & {} & {}\\
e^\pm_{m+n,1}(z) & \cdots & e^\pm_{m+n,m+n-1}(z) & 1
\end{array}
\right )
\left (
\begin{array}{ccc}
k^\pm_1(z) & \cdots & 0\\
\vdots & \ddots & \vdots\\
0 & \cdots & k^\pm_{m+n}(z)
\end{array}
\right )\no\\
& & \times \left (
\begin{array}{ccccc}
1 & f^\pm_{1,2}(z) & f^\pm_{1,3}(z) & \cdots & f^\pm_{1,m+n}(z)\\
\vdots & \ddots & \cdots & {} & \vdots\\
{} & {} & {} & {} & f^\pm_{m+n-1,m+n}(z)\\
0 & {} & {} & {} & 1
\end{array}
\right ),
\eea
where $e^\pm_{i,j}(z),~f^\pm_{j,i}(z)$ and $k^\pm_i(z) ~(i>j)$ are 
elements in the super RS algebra and $k^\pm_i(z)$ are invertible. 
Let
\bea
X^-_i(z)&=&f^+_{i,i+1}(z_+)-f^-_{i,i+1}(z_-),\no\\
X^+_i(z)&=&e^+_{i+1,i}(z_-)-e^-_{i+1,i}(z_+),
\eea
where $z_\pm=zq^{\pm{c\over 2}}$, then $q^{\pm{c\over 2}},\;
X^\pm_i(z),\;k^\pm_j(z),~i=1,2,\cdots,m+n-1,\;j=1,2,\cdots,m+n$
give the defining relations of quantum affine superalgebra $\uqgh$.
\end{Theorem}

The Gauss decomposition implies that the elements $e^\pm_{i,j}(z),\;
f^\pm_{j,i}(z) ~(i>j)$ and $k^\pm_i(z)$ are uniquely determined by
$L^\pm(z)$. In the following we will denote $f^\pm_{i,i+1}(z),\;
e^\pm_{i+1,i}(z)$ as $f^\pm_i(z),\;e^\pm_i(z)$, respectively.

The following matrix equations can be deduced from (\ref{RLL-LLR1}):
\bea
R_{21}({z\over w})\t L^\pm_2(z)\t L_1^\pm(w)&=&
  L^\pm_1(w)\t L^\pm_2(z)\t R_{21}({z\over w}),\label{RLL-LLR2}\\
R_{21}({z_-\over w_+})\t L^-_2(z)\t L_1^+(w)&=&
  L^+_1(w)\t L^-_2(z)\t R_{21}({z_+\over w_-}),\label{RLL-LLR3}\\
\t L^\pm_2(z)^{-1}\t L^\pm_1(w)^{-1}R_{21}({z\over w})&=&
  R_{21}({z\over w})L^\pm_1(w)^{-1}\t L^\pm_2(z)^{-1}
  \t,\label{RLL-LLR4}\\
\t L^+_2(z)^{-1}\t L^-_1(w)^{-1}R_{21}({z_+\over w_-})&=&
  R_{21}({z_-\over w_+})L^-_1(w)^{-1}\t L^+_2(z)^{-1}
  \t,\label{RLL-LLR5}\\
L^\pm_1(w)^{-1}R_{21}({z\over w})\t L^\pm_2(z)\t&=&
  \t L^\pm_2(z)\t R_{21}({z\over
  w})L^\pm_1(w)^{-1},\label{RLL-LLR6}\\
L^-_1(w)^{-1}R_{21}({z_+\over w_-})\t L^+_2(z)\t&=&
  \t L^+_2(z)\t R_{21}({z_-\over
  w_+})L^-_1(w)^{-1},\label{RLL-LLR7}\\
L^+_1(w)^{-1}R_{21}({z_-\over w_+})\t L^-_2(z)\t&=&
  \t L^-_2(z)\t R_{21}({z_+\over
  w_-})L^+_1(w)^{-1},\label{RLL-LLR8}
\eea
where $R_{21}({z\over w})=R({w\over z})^{-1}$. As in (\ref{RLL-LLR1}), the
multiplications in (\ref{RLL-LLR2} -- \ref{RLL-LLR8}) are usual matrix
multiplications.

\sect{The $m=1,\;n=1$ Case: $U_q[gl(1|1)^{(1)}]$}

For the simplest supersymmetric case $U_q[gl(1|1)^{(1)}]$, 
$L^\pm(z)$ take the forms 
\beq
L^\pm(z)=\left (
\begin{array}{cc}
k^\pm_1(z) & k^\pm_1(z)f^\pm_1(z)\\
e^\pm_1(z)k^\pm_1(z) & k^\pm_2(z)+e^\pm_1(z)k^\pm_1(z)f^\pm_1(z)
\end{array}
\right ).\label{L-gl11}
\eeq
R-matrix and $\t$ are given by,
\beq
R(\frac{z}{w})=\lt(
\begin{array}{cccc}
1 & 0 & 0 & 0\\
0 & \frac{z-w}{zq-wq^{-1}} & \frac{z(q-q^{-1})}{zq-wq^{-1}} & 0\\
0 & \frac{w(q-q^{-1})}{zq-wq^{-1}} & \frac{z-w}{zq-wq^{-1}} & 0\\
0 & 0 & 0 & -\frac{wq-zq^{-1}}{zq-wq^{-1}}
\end{array}
\rt)
\eeq
and $\t=diag(1,1,1,-1)$, respectively.

Using (\ref{theta}), (\ref{r-matrix}) and (\ref{RLL-LLR1},
\ref{RLL-LLR2} -- \ref{RLL-LLR8}), and by 
similar calculations as to the non-super case \cite{Din93}, we obtain
a super RS algebra, which leads to, by means of theorem \ref{df},
\begin{Definition}\label{simplest}:
$U_q[gl(1|1)^{(1)}]$ is an associative algebra with unit 1 and the
Drinfeld current generators: $X^\pm_1(z),~k^\pm_i(z)~(i=1,2)$, a central
element $c$ and a nonzero complex parameter $q$. $k^\pm_i(z)$ are 
invertible. The gradings of the generators are:
$[X^\pm_1(z)]=1$ and $[k^\pm_i(z)]=0=[c]$. The relations read,
\bea
k^\pm_i(z)k^\pm_j(w)&=&k^\pm_j(w)k^\pm_i(z),~~i,\;j=1,2,\no\\
k^+_1(z)k^-_1(w)&=&k^-_1(w)k^+_1(z),\no\\
{{w_-q-z_+q^{-1}}\over{z_+q-w_-q^{-1}}}k^+_2(z)k^-_2(w)&=&
  {{w_+q-z_-q^{-1}}\over{z_-q-w_+q^{-1}}}k^-_2(w)k^+_2(z),\no\\
{{z_\pm-w_\mp}\over{z_\pm q-w_\mp q^{-1}}}k^\mp_2(w)^{-1}k^\pm_1(z)&=&
  {{z_\mp-w_\pm}\over{z_\mp q-w_\pm
  q^{-1}}}k^\pm_1(z)k^\mp_2(w)^{-1},\no\\
k^\pm_i(z)^{-1}X^-_1(w)k^\pm_i(z)&=&
  \frac{z_\mp q-wq^{-1}}{z_\mp-w}X^-_1(w),\no\\
k^\pm_i(z)X^+_1(w)k^\pm_i(z)^{-1}&=&
  \frac{z_\pm q-wq^{-1}}{z_\pm-w}X^+_1(w),\no\\
\{X^\pm_1(z),X^\pm_1(w)\}&=&0,\no\\
\{X^+_1(z),X^-_1(w)\}&=&(q-q^{-1})\lt(\d({w\over z}q^{c})k^+_2(w_+)
   k^+_1(w_+)^{-1}\rt.\no\\
& &~~~~~\lt. -\d({w\over z}q^{-c})k^-_2(z_+)k^-_1(z_+)^{-1}
   \rt),\label{gl11-rs}
\eea
where $\{X,Y\}\equiv XY+YX$ denotes an anti-commutator and
\beq
\d(z)=\sum_{k\in {\bf Z}}\,z^k
\eeq
is a formal series.
\end{Definition}

\begin{Theorem}\label{hopf-gl11}: The algebra $U_q[gl(1|1)^{(1)}]$
defined by (\ref{gl11-rs}) has a Hopf algebra structure, which is given
by the following formulae.\\
{\bf Coproduct} $\D$
\bea
&&\D(q^c)=q^c\otimes q^c,\no\\
&&\D(k^+_i(z))=k^+_i(zq^{{c_2\over 2}})\otimes k^+_i(zq^{-{c_1\over
   2}}),\no\\ 
&&\D(k^-_i(z))=k^-_i(zq^{-{c_2\over 2}})\otimes k^-_i(zq^{{c_1\over
   2}}),\no\\ 
&&\D(X^+_1(z))=X^+_1(z)\otimes 1+\psi_1(zq^{{c_1\over 2}})\otimes
   X^+_1(zq^{c_1}),\no\\
&&\D(X^-_1(z))=1\otimes X^-_1(z)+X^-_1(zq^{c_2})\otimes
  \phi_1(zq^{{c_2\over 2}}), \label{coproduct-gl11}
\eea
where $c_1=c\otimes 1$, $c_2=1\otimes c$,
$\psi_1(z)=k^-_2(z)k^-_1(z)^{-1}$ and
$\phi_1(z)=k^+_2(z)k^+_1(z)^{-1}$.\\
{\bf Counit} $\e$
\beq
\e(q^c)=1,~~~~\e(k^\pm_i(z))=1,~~~~\e(X^\pm_1(z))=0.\label{counit-gl11}
\eeq
{\bf Antipode} $S$
\bea
&&S(q^c)=q^{-c},~~~~S(k^\pm_i(z))=k^\pm_i(z)^{-1},\;i=1,2,\no\\
&&S(X^+_1(z))=-\psi_1(zq^{-{c\over 2}})^{-1}X^+_1(zq^{-c}),\no\\
&&S(X^-_1(z))=-X^-_1(zq^{-c})\phi_1(zq^{-{c\over
    2}})^{-1}.\label{antipode-gl11}
\eea
\end{Theorem}

\noindent {\bf Proof}. The proof is rather elementary. We nevertheless
present the details since the proof for the general case in next section
is similar. Care has to be taken of the gradings in
tensor product multiplications and also in extending the antipode to
the whole algebra.
\bea
&&\D\lt(\frac{w_-q-z_+q^{-1}}{z_+q-w_-q^{-1}}k^+_2(z)k^-_2(w)\rt)=
  \frac{wq^{\frac{-c_1-c_2}{2}+1}-zq^{\frac{c_1+c_2}{2}-1}}
  {zq^{\frac{c_1+c_2}{2}+1}-wq^{\frac{-c_1-c_2}{2}-1}}\no\\
&&~~~~~~~~\times \lt(k^+_2(zq^{\frac{c_2}{2}})\otimes k^+_2
  (zq^{-\frac{c_1}{2}})\rt)
  \lt(k^-_2(wq^{-\frac{c_2}{2}})\otimes k^-_2(wq^{\frac{c_1}{2}})
  \rt)\no\\
&&~~~~~~~~=\frac{wq^{\frac{c_1+c_2}{2}+1}-zq^{\frac{-c_1-c_2}{2}-1}}
  {zq^{\frac{-c_1-c_2}{2}+1}-wq^{\frac{c_1+c_2}{2}-1}}\no\\
&&~~~~~~~~\times \lt(k^-_2(wq^{-\frac{c_2}{2}})\otimes k^-_2
  (wq^{\frac{c_1}{2}})\rt)
  \lt(k^+_2(zq^{\frac{c_2}{2}})\otimes k^+_2(zq^{-\frac{c_1}{2}})\rt)\no\\
&&~~~~~~~~=\D\lt(\frac{w_+q-z_-q^{-1}}{z_-q-w_+q^{-1}}k^-_2(w)k^+_2(z)\rt),
\eea
where the relation between $k^+_2(z)$ and $k^-_2(w)$ has been used.
\bea
&&\D\lt(\frac{z_+-w_-}{z_+q-w_-q^{-1}}k^+_1(z)k^-_2(w)\rt)=
  \frac{zq^{\frac{c_1+c_2}{2}}-wq^{\frac{-c_1-c_2}{2}}}
  {zq^{\frac{c_1+c_2}{2}+1}-wq^{\frac{-c_1-c_2}{2}-1}}\no\\
&&~~~~~~~~\times \lt(k^+_1(zq^{\frac{c_2}{2}})
  \otimes k^+_1(zq^{-\frac{c_1}{2}})\rt)
  \lt(k^-_2(wq^{-\frac{c_2}{2}})\otimes k^-_2(wq^{\frac{c_1}{2}})
  \rt)\no\\
&&~~~~~~~~ = \frac{zq^{\frac{-c_1-c_2}{2}}-wq^{\frac{c_1+c_2}{2}}}
  {zq^{\frac{-c_1-c_2}{2}+1}-wq^{\frac{c_1+c_2}{2}-1}}\no\\
&&~~~~~~~~ \times \lt(k^-_2(wq^{-\frac{c_2}{2}})
  \otimes k^-_2(wq^{\frac{c_1}{2}})\rt)
  \lt(k^+_1(zq^{\frac{c_2}{2}})\otimes k^+_1(zq^{-\frac{c_1}{2}})\rt)\no\\
&&~~~~~~~~ = \D\lt(\frac{z_--w_+}{z_-q-w_+q^{-1}}k^-_2(w)k^+_1(z)\rt).
\eea
The relation between $\D(k^-_1(z))$ and $\D(k^+_2(w))$ can be proved 
along a simliar line.

The fifth and sixth equations in (\ref{gl11-rs}) are equivalent to the
relations,
\bea
\phi_1(z)X^\pm_1(w)\phi_1(z)^{-1}=X^\pm_1(w),\no\\
\psi_1(z)X^\pm_1(w)\psi_1(z)^{-1}=X^\pm_1(w)
\eea
which are easily seen to be preserved by the coproduct, thanks
to the commutativity of $\phi_1(z)$ and $\psi_1(w)$.
\bea
\D\lt(\{X^+_1(z),X^+_1(w)\}\rt)&=&\D(X^+_1(z))\D(X^+_1(w))+\D(X^+_1(w))
    \D(X^+_1(z))\no\\
&=&X^+_1(z)X^+_1(w)\otimes 1+X^+_1(z)\psi_1(wq^{\frac{c_1}{2}})
    \otimes X^+_1(wq^{c_1})\no\\
& &    -\psi_1(zq^{\frac{c_1}{2}})X^+_1(w)\otimes
    X^+_1(zq^{c_1})\no\\
& &+X^+_1(w)X^+_1(z)\otimes 1+X^+_1(w)\psi_1(zq^{\frac{c_1}{2}})
    \otimes X^+_1(zq^{c_1})\no\\
& &    -\psi_1(wq^{\frac{c_1}{2}})X^+_1(z)
    \otimes X^+_1(wq^{c_1})\no\\
& &+\psi_1(zq^{\frac{c_1}{2}})\psi_1(wq^{\frac{c_1}{2}})\otimes
    X^+_1(zq^{c_1})X^+_1(wq^{c_1})\no\\
& &    +\psi_1(wq^{\frac{c_1}{2}})
    \psi_1(zq^{\frac{c_1}{2}})\otimes
    X^+_1(wq^{c_1})X^+_1(zq^{c_1})\no\\
&=&\{X^+_1(z),X^+_1(w)\}\otimes 1+\psi_1(zq^{\frac{c_1}{2}})
    \psi_1(wq^{\frac{c_1}{2}})\no\\
& &    \otimes \{X^+_1(zq^{c_1}),
    X^+_1(wq^{c_1})\}=0.
\eea
$\D(\{X^-_1(z),X^-_1(w)\})=0$ can be proved similarly.
\bea
\D\lt(\{X^+_1(z),X^-_1(w)\}\rt)&=&\D(X^+_1(z))\D(X^-_1(w))+\D(X^-_1(w))
    \D(X^+_1(z))\no\\
&=&X^+_1(z)\otimes X^-_1(w)+X^+_1(z)X^-_1(wq^{c_2})\otimes
    \phi_1(wq^{\frac{c_2}{2}})\no\\
& &    +\psi_1(zq^{\frac{c_1}{2}})\otimes
    X^+_1(zq^{c_1})X^-_1(w)\no\\
& &-X^+_1(z)\otimes X^-_1(w)+\psi_1(zq^{\frac{c_1}{2}})\otimes
    X^-_1(w)X^+_1(zq^{c_1})\no\\
& &    +X^-_(wq^{c_2})X^+_1(z)\otimes
    \phi_1(wq^{\frac{c_2}{2}})\no\\
& &-\psi_1(zq^{\frac{c_1}{2}})X^-_1(wq^{c_2})\otimes
    X^+_1(zq^{c_1})\phi_1(wq^{\frac{c_2}{2}})\no\\
& &    +X^-_1(wq^{c_2})
    \psi_1(zq^{\frac{c_1}{2}})X^+_1(zq^{c_1})\no\\
&=&[X^+_1(z),X^-_1(wq^{c_2})]\otimes\phi_1(wq^{c_2})
    +\psi_1(zq^{\frac{c_1}{2}})\otimes [X^+_1(zq^{c_1}),X^-_1(w)]\no\\
&=&(q-q^{-1})\lt(\d(\frac{w}{z}q^{c_1+c_2})\phi_1(wq^{\frac{c_2}{2}+
   \frac{c_1+c_2}{2}})\phi_1(wq^{-\frac{c_1}{2}+\frac{c_1+c_2}{2}})\rt.\no\\
& &\lt. -\d(\frac{w}{z}q^{-c_1-c_2})\psi_1(zq^{\frac{-c_2}{2}+
   \frac{c_1+c_2}{2}})\psi_1(wq^{\frac{c_1}{2}+\frac{c_1+c_2}{2}})\rt)\no\\
&=&(q-q^{-1})\D\lt(\d(\frac{w}{z}q^c)\phi_1(w_+)
   -\d(\frac{w}{z}q^{-c})\psi_1(z_+)\rt).
\eea
We have therefore proved that the comultiplication is an algebra 
homomorphism.
\bea
S\lt(\{X^+_1(z), X^-_1(w)\}\rt)&=&-S(X^-_1(w))S(X^+_1(z))-S(X^+_1(z))
     S(X^-_1(w))\no\\
&=&-\psi_1(zq^{-\frac{c}{2}})^{-1}\phi_1(wq^{-\frac{c}{2}})^{-1}
    [X^+_1(zq^{-c}), X^-_1(wq^{-c})]\no\\
&=&-\psi_1(zq^{-\frac{c}{2}})^{-1}\phi_1(wq^{-\frac{c}{2}})^{-1}
   (q-q^{-1})\no\\
& &\times   \lt(\d(\frac{w}{z}q^c)\phi_1(w^{-\frac{c}{2}})
   -\d(\frac{w}{z}q^{-c})\psi_1(zq^{-\frac{c}{2}})\rt)\no\\
&=&(q-q^{-1})\lt(\d(\frac{w}{z}q^{-c})\phi_1(w^{-\frac{c}{2}})^{-1}
   -\d(\frac{w}{z}q^{c})\psi_1(zq^{-\frac{c}{2}})^{-1}\rt)\no\\
&=&(q-q^{-1})S\lt(\d(\frac{w}{z}q^{c})\phi_1(w^{\frac{c}{2}})
   -\d(\frac{w}{z}q^{-c})\psi_1(zq^{\frac{c}{2}})\rt).
\eea
We can prove in the same manner that other relations are also preserved
by the antipode.

Let $M:~U_q[gl(1|1)^{(1)}]\otimes U_q[gl(1|1)^{(1)}]\rightarrow
U_q[gl(1|1)^{(1)}]$ be a multiplication. Then we can easily check
\bea
&&M(1\otimes\e)\D=id=M(\e\otimes 1)\D,\no\\
&&M(1\otimes S)\D=\e=M(S\otimes 1)\D.
\eea
Thus we have shown that the coproduct, the counit and the antipode
give a Hopf algebra structure.

\sect{General Case: $U_q[gl(m|n)^{(1)}]$}

The generalization to the general case $U_q[gl(m|n)^{(1)}]$ is
more or less straightforward
by using (\ref{r-matrix}), (\ref{RLL-LLR1}, 
\ref{RLL-LLR2} -- \ref{RLL-LLR8}) and theorem \ref{df}. 
Similar to the non-supersymmetric cases \cite{Din93,Fan97}, 
this is achieved by induction on
$m$ and $n$. Tedious but direct computations give rise to

\begin{Definition}\label{general}: $U_q[gl(m|n)^{(1)}]$ is an
associative algebra with unit 1 and Drinfeld current generators: $X^\pm_i(z),~
k^\pm_j(z),~i=1,2,\cdots,m+n-1,~j=1,2,\cdots,m+n$, a central element
$c$ and a nonzero complex parameter $q$. $k^\pm_i(z)$ are invertible.
The grading of the generators are: $[X^\pm_m(z)]=1$ and zero otherwise.
The defining relations are given by
\bea
k^\pm_i(z)k^\pm_j(w)&=&k^\pm_j(w)k^\pm_i(z),~~i\neq j\no\\
k^+_i(z)k^-_i(w)&=&k^-_i(w)k^+_i(z),~~i\leq m,\no\\
{{w_-q-z_+q^{-1}}\over{z_+q-w_-q^{-1}}}k^+_i(z)k^-_i(w)&=&
  {{w_+q-z_-q^{-1}}\over{z_-q-w_+q^{-1}}}k^-_i(w)k^+_i(z),~~m<i\leq
  m+n,\no\\
{{z_\pm-w_\mp}\over{z_\pm q-w_\mp q^{-1}}}k^\mp_i(w)^{-1}k^\pm_j(z)&=&
  {{z_\mp-w_\pm}\over{z_\mp q-w_\pm
  q^{-1}}}k^\pm_j(z)k^\mp_i(w)^{-1},~~i>j,\no\\
k^\pm_j(z)^{-1}X^-_i(w)k^\pm_j(z)&=&X^-_i(w),~~j-i\leq -1,\no\\
k^\pm_j(z)^{-1}X^+_i(w)k^\pm_j(z)&=&X^+_i(w),~~j-i\leq -1,~~{\rm or}\no\\
k^\pm_j(z)^{-1}X^-_i(w)k^\pm_j(z)&=&X^-_i(w),~~j-i\geq 2,\no\\
k^\pm_j(z)^{-1}X^+_i(w)k^\pm_j(z)&=&X^+_i(w),~~j-i\geq 2,\no\\
k^\pm_i(z)^{-1}X^-_i(w)k^\pm_i(z)&=&
  \frac{z_\mp q-wq^{-1}}{z_\mp-w}X^-_i(w),~~i<m,\no\\
k^\pm_i(z)^{-1}X^-_i(w)k^\pm_i(z)&=&
  \frac{z_\mp q^{-1}-wq}{z_\mp-w}X^-_i(w),~~m<i\leq m+n-1,\no\\
k^\pm_{i+1}(z)^{-1}X^-_i(w)k^\pm_{i+1}(z)&=&
  \frac{z_\mp q^{-1}-wq}{z_\mp-w}X^-_i(w),~~i<m,\no\\
k^\pm_{i+1}(z)^{-1}X^-_i(w)k^\pm_{i+1}(z)&=&
  \frac{z_\mp q-wq^{-1}}{z_\mp-w}X^-_i(w),~~m<i\leq m+n-1,\no\\
k^\pm_i(z)X^+_i(w)k^\pm_i(z)^{-1}&=&
  \frac{z_\pm q-wq^{-1}}{z_\pm-w}X^+_i(w),~~i<m,\no\\
k^\pm_i(z)X^+_i(w)k^\pm_i(z)^{-1}&=&
  \frac{z_\pm q^{-1}-wq}{z_\pm-w}X^+_i(w),~~m<i\leq m+n-1,\no\\
k^\pm_{i+1}(z)X^+_i(w)k^\pm_{i+1}(z)^{-1}&=&
  \frac{z_\pm q^{-1}-wq}{z_\pm-w}X^+_i(w),~~i<m,\no\\
k^\pm_{i+1}(z)X^+_i(w)k^\pm_{i+1}(z)^{-1}&=&
  \frac{z_\pm q-wq^{-1}}{z_\pm-w}X^+_i(w),~~m<i\leq m+n-1,\no\\
k^\pm_i(z)^{-1}X^-_m(w)k^\pm_i(z)&=&
  \frac{z_\mp q-wq^{-1}}{z_\mp-w}X^-_m(w),~~i=m,\;m+1,\no\\
k^\pm_i(z)X^+_m(w)k^\pm_i(z)^{-1}&=&
  \frac{z_\pm q-wq^{-1}}{z_\pm-w}X^+_m(w),~~i=m,\;m+1,\no\\
(zq^{\mp 1}-wq^{\pm 1})X^\mp_i(z)X^\mp_i(w)&=&(zq^{\pm 1}-wq^{\mp 1})
  X^\mp_i(w)X^\mp_i(z),~~i<m,\no\\
(wq^{\mp 1}-zq^{\pm 1})X^\mp_i(z)X^\mp_i(w)&=&(wq^{\pm 1}-zq^{\mp 1})
  X^\mp_i(w)X^\mp_i(z),~~m<i\leq m+n-1,\no\\
\{X^\pm_m(z),X^\pm_m(w)\}&=&0,\no\\
(z-w)X^+_i(z)X^+_{i+1}(w)&=&(zq-wq^{-1})
  X^+_{i+1}(w)X^+_i(z),~~i<m,\no\\
(w-z)X^+_i(z)X^+_{i+1}(w)&=&(wq-zq^{-1})
  X^+_{i+1}(w)X^+_i(z),~~m\leq i\leq m+n-1,\no\\
(zq-wq^{-1})X^-_i(z)X^-_{i+1}(w)&=&(z-w)
  X^-_{i+1}(w)X^-_i(z),~~i<m,\no\\
(wq-zq^{-1})X^-_i(z)X^-_{i+1}(w)&=&(w-z)
  X^-_{i+1}(w)X^-_i(z),~~m\leq i\leq m+n-1,\no\\
{[X^+_i(z),X^-_j(w)]}&=&-(q-q^{-1})\d_{ij}\lt(
   \d({w\over z}q^{c})k^+_{i+1}(w_+)
   k^+_i(w_+)^{-1}\rt.\no\\
& &~~~~~\lt. -\d({w\over z}q^{-c})k^-_{i+1}(z_+)k^-_i(z_+)^{-1}
  \rt),~~i,j\neq m,\no\\
\{X^+_m(z),X^-_m(w)\}&=&(q-q^{-1})\lt(
   \d({w\over z}q^{c})k^+_{m+1}(w_+)
   k^+_m(w_+)^{-1}\rt.\no\\
& &~~\lt. -\d({w\over z}q^{-c})k^-_{m+1}(z_+)k^-_m(z_+)^{-1}\rt),
  \label{glmn-rs}
\eea
where $[X,Y]\equiv XY-YX$ stands for a commutator and $\{X,Y\}\equiv XY+YX$
for an anti-commutator, together with
the following Serre and extra Serre \cite{Sch92,Yam96} relations:
\bea
&&\{X^\pm_i(z_1)X^\pm_i(z_2)X^\pm_{i+1}(w)-(q+q^{-1})X^\pm_i(z_1)
  X^\pm_{i+1}(w)X^\pm_i(z_2)\no\\
&&~~~~~~~~~~  +X^\pm_{i+1}(w)X^\pm_i(z_1)X^\pm
  _i(z_2)\}+\{z_1\leftrightarrow z_2\}=0,~~i\neq m,
  \label{serre1}\\
&&\{X^\pm_{i+1}(z_1)X^\pm_{i+1}(z_2)X^\pm_{i}(w)-(q+q^{-1})X^\pm_{i+1}(z_1)
  X^\pm_{i}(w)X^\pm_{i+1}(z_2)\no\\
&&~~~~~~~~~~  +X^\pm_{i}(w)X^\pm_{i+1}(z_1)X^\pm
  _{i+1}(z_2)\}+\{z_1\leftrightarrow z_2\}=0,~~i\neq m-1,
  \label{serre2}\\
&&\{(z_1q^{\mp 1}-z_2q^{\pm 1})[X^\pm_m(z_1)X^\pm_m(z_2)X^\pm_{m-1}(w)-
  (q+q^{-1})X^\pm_m(z_1)X^\pm_{m-1}(w)X^\pm_m(z_2)\no\\
&&~~~~~~~~~~  +X^\pm_{m-1}(w)X^\pm_m(z_1)X^\pm
    _m(z_2)]\}+\{z_1\leftrightarrow z_2\}=0,
    \label{serre3}\\
&&\{(z_2q^{\mp 1}-z_1q^{\pm 1})[X^\pm_m(z_1)X^\pm_m(z_2)X^\pm_{m+1}(w)-
  (q+q^{-1})X^\pm_m(z_1)X^\pm_{m+1}(w)X^\pm_m(z_2)\no\\
&&~~~~~~~~~~  +X^\pm_{m+1}(w)X^\pm_m(z_1)X^\pm
     _m(z_2)]\}+\{z_1\leftrightarrow z_2\}=0,
     \label{serre4}\\
&&\{(z_1q^{\mp 1}-z_2q^{\pm 1})[X^\pm_m(z_1)X^\pm_m(z_2)X^\pm_{m-1}(w_1)
  X^\pm_{m+1}(w_2)\no\\
&&~~~~~~~~~~  -(q+q^{-1})X^\pm_m(z_1)X^\pm_{m-1}(w_1)X^\pm_m(z_2)
  X^\pm_{m+1}(w_2)]\no\\
&&~~~~~~~~~~+(z_1+z_2)(q^{\mp 1}-q^{\pm 1})X^\pm_{m-1}(w_1)X^\pm_m(z_1)X^\pm
            _m(z_2)X^\pm_{m+1}(w_2)\no\\
&&~~~~~~~~~~+(z_2q^{\mp 1}-z_1q^{\pm 1})[-(q+q^{-1})X^\pm_{m-1}(w_1)
  X^\pm_m(z_1)X^\pm_{m+1}(w_2)X^\pm_m(z_2)\no\\
&&~~~~~~~~~~+X^\pm_{m-1}(w_1)X^\pm_{m+1}(w_2)X^\pm_m(z_1)
  X^\pm_m(z_2)]\}+\{z_1\leftrightarrow z_2\}=0.
  \label{extra-serre}
\eea
\end{Definition}

\begin{Theorem}\label{hopf-glmn}: The algebra $U_q[gl(m|n)^{(1)}]$
given by definition \ref{general} has a Hopf algebra structure, which
is given by the following formulae.\\
{\bf Coproduct} $\D$
\bea
&&\D(q^c)=q^c\otimes q^c,\no\\
&&\D(k^+_j(z))=k^+_j(zq^{{c_2\over 2}})\otimes k^+_j(zq^{-{c_1\over
   2}}),\no\\ 
&&\D(k^-_j(z))=k^-_j(zq^{-{c_2\over 2}})\otimes k^-_j(zq^{{c_1\over
   2}}),~~j=1,2,\cdots,m+n\no\\ 
&&\D(X^+_i(z))=X^+_i(z)\otimes 1+\psi_i(zq^{{c_1\over 2}})\otimes
   X^+_i(zq^{c_1}),\no\\
&&\D(X^-_i(z))=1\otimes X^-_i(z)+ X^-_i(zq^{c_2})\otimes
  \phi_i(zq^{{c_2\over 2}}),~~i=1,2,\cdots,m+n-1, \label{coproduct-glmn}
\eea
where $c_1=c\otimes 1,~c_2=1\otimes c,~\psi_i(z)=k^-_{i+1}(z)k^-_i(z)^{-1}$
and $\phi_i(z)=k^+_{i+1}(z)k^+_i(z)^{-1}$.\\
{\bf Counit} $\e$
\beq
\e(q^c)=1,~~~~\e(k^\pm_j(z))=1,~~~~\e(X^\pm_i(z))=0.\label{counit-glmn}
\eeq
{\bf Antipode} $S$
\bea
&&S(q^c)=q^{-c},~~~~S(k^\pm_j(z))=k^\pm_j(z)^{-1},\no\\
&&S(X^+_i(z))=-\psi_i(zq^{-{c\over 2}})^{-1}X^+_i(zq^{-c}),\no\\
&&S(X^-_i(z))=-X^-_i(zq^{-c})\phi_i(zq^{-{c\over
    2}})^{-1}.\label{antipode-glmn}
\eea
\end{Theorem}

\noindent{\bf Proof}: Similar to the case in the previous
section for $U_q[gl(1|1)^{(1)}]$, this theorem is proved by direct
calculations.

\vskip.3in
\noindent {\bf Acknowledgements.} This work is supported by
Australian Research Council, and in part by University of Queensland
New Staff Research Grant and External Support Enabling Grant.


\end{document}